\documentclass[reqno]{amsart}
\usepackage{amsrefs}
\usepackage{amssymb}
\usepackage{graphicx}
\usepackage{url}
\usepackage{array}
\usepackage{booktabs,multirow,tabularx}

\usepackage{graphicx}
\usepackage{pgf,tikz,circuitikz,subfig}
\usetikzlibrary{arrows,shapes,positioning}

\newcolumntype{L}[1]{>{\raggedright\let\newline\\\arraybackslash\hspace{0pt}}m{#1}}

\usepackage{xcolor}

\newtheorem{thm}{Theorem}

\DeclareMathOperator{\subj}{s.t.}

\newcommand{\ceiling}[1]{\left\lceil#1\right\rceil}
\newcommand{\floor}[1]{\left\lfloor#1\right\rfloor}

\newcommand\TS{\rule{0pt}{2.6ex}} 
\newcommand\BS{\rule[-1.2ex]{0pt}{0pt}} 

\allowdisplaybreaks

\begin{document}

\title{Small polygons with large area}

\author[C. Bingane]{Christian Bingane}
\address{Department of Mathematics and Industrial Engineering, Polytechnique Montreal, Montreal, QC, Canada}
\email{christian.bingane@polymtl.ca}

\author[M.~J. Mossinghoff]{Michael J. Mossinghoff}
\address{Center for Communications Research, Princeton, NJ, USA}
\email{m.mossinghoff@idaccr.org}

\date\today

\begin{abstract}
A polygon is \textit{small} if it has unit diameter.
The maximal area of a small polygon with a fixed number of sides $n$ is not known when $n$ is even and $n\geq14$.
We determine an improved lower bound for the maximal area of a small $n$-gon for this case.
The improvement affects the $1/n^3$ term of an asymptotic expansion; prior advances affected less significant terms.
This bound cannot be improved by more than $O(1/n^3)$.
For $n=6$, $8$, $10$, and $12$, the polygon we construct has maximal area.
\end{abstract}

\maketitle

\section{Introduction}\label{sectionIntroduction}

A polygon is said to be \textit{small} if it has diameter~$1$.
Reinhardt \cite{Reinhardt22} first studied two extremal problems for small polygons a century ago:
determining the maximal area for a small polygon with $n$ sides, and determining the maximal perimeter for a small convex polygon with $n$ sides.
These are sometimes referred to as \textit{isodiametric problems} for polygons in the literature.
See \cite{AHM07,AHM09} for a survey of work on these problems and related ones.
We focus on the area problem in this paper.

Reinhardt proved that the regular small polygon alone has maximal area when $n$ is odd, and that this polygon is never optimal when $n$ is even and $n\ge 6$.
It is straightforward to show that there are infinitely many different small quadrilaterals with maximal area, including the square, and the optimal hexagon was first determined by Graham in 1975 \cite{Graham}.
The optimal octagon was established by Audet et al.\ in 2002 \cite{AHMX}, and the cases $n=10$ and $n=12$ were resolved by Henrion and Messine in 2013 \cite{HenrionMessine}.
The problem remains open for larger $n$.

In 2006, Foster and Szabo \cite{FosterSzabo} proved that the area $A(P_n)$ of a small polygon $P_n$ having an even number of sides $n$ satisfies $A(P_n) < \overline{A}_n$, where
\begin{equation}\label{eqnAreaBound}
\begin{aligned}
\overline{A}_n &= \frac{n}{2}\sin\left(\frac{\pi}{n}\right) - \frac{n-1}{2}\tan\left(\frac{\pi}{2n-2}\right)\\
&= \frac{\pi}{4} - \frac{5\pi^3}{48n^2} - \frac{\pi^3}{24n^3} + O\left(\frac{1}{n^4}\right).
\end{aligned}
\end{equation}
Since the area of the small regular polygon $R_n$ with $n$ even is given by $\frac{n}{8}\sin(2\pi/n)$, it follows easily that
\[
\overline{A}_n - A(R_n) = \frac{\pi^3}{16n^2} + O\left(\frac{1}{n^3}\right).
\]
In 2005, the second author \cite{Mossinghoff05} described a construction for a polygon $M_n$ with an even number of sides $n$ for which
\[
\overline{A}_n - A(M_n) = \frac{a_3\pi^3}{n^3} + O\left(\frac{1}{n^4}\right),
\]
with
\[
a_3 = \frac{5303-456\sqrt{114}}{5808} = 0.07476796\ldots\,.
\]
The first author \cite{BinganeA} recently improved this, constructing a polygon $B_n$ for each even $n\geq6$ satisfying
\[
A(B_n) - A(M_n) = \frac{a_5\pi^3}{n^5}  + O\left(\frac{1}{n^6}\right),
\]
with $a_5=0.25097\ldots$ when $n\equiv2\bmod4$ and $a_5=0.35411\ldots$ when $n\equiv0\bmod4$.
In this paper, we generalize the latter construction to produce small polygons $Q_n$ that exhibit an improvement in the $1/n^3$ term for the area problem.
We establish the following result.

\begin{thm}\label{thmArea}
Let $n\geq6$ be an even integer, let $B_n$ denote the small $n$-gon from \cite{BinganeA}, and let $\overline{A}_n$ denote the upper bound on the area of a small $n$-gon given by \eqref{eqnAreaBound}.
There exists a small $n$-gon $Q_n$ satisfying
\[
\overline{A}_n - A(Q_n) = \frac{\delta\pi^3}{n^3} + O\left(\frac{1}{n^4}\right) < \frac{8\pi^3}{109n^3} + O\left(\frac{1}{n^4}\right),
\]
with $\delta=0.0733883168\ldots$\,, so
\[
A(Q_n) - A(B_n) > \frac{\pi^3}{725n^3} +  O\left(\frac{1}{n^4}\right).  
\]
Moreover, $Q_n$ is the optimal small polygon for $n\leq12$.
\end{thm}

This article is organized in the following way.
Section~\ref{secPrior} establishes some notation and describes some prior constructions.
Section~\ref{secProof} describes the new construction and proves Theorem~\ref{thmArea}.
Section~\ref{secSmallN} reports on the computation of optimal small $n$-gons for a number of even $n$ assuming of an axis of symmetry, and compares the results of our construction with these polygons.

\section{Prior constructions}\label{secPrior}

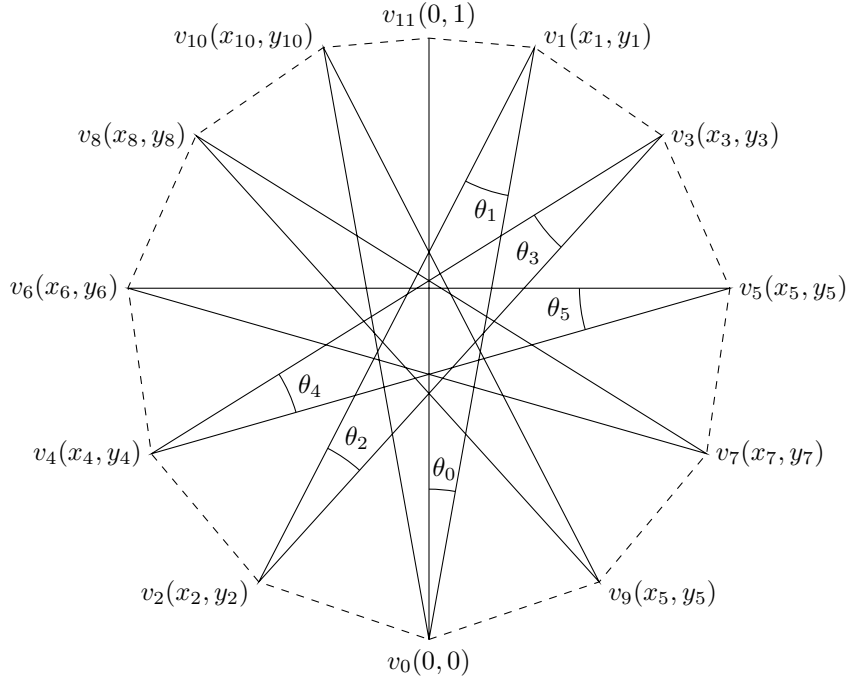
\begin{figure}
\caption{Polygon and skeleton for $n=12$.}\label{figSkeleton}
	\centering
	\begin{tikzpicture}[scale=8]
		\draw[dashed] (0,0) node[below]{$v_0(0,0)$} -- (0.2828,0.0957) node[right,yshift=-1ex]{$v_9(x_5,y_5)$} -- (0.4614,0.3088) node[right]{$v_7(x_7,y_7)$} -- (0.5000,0.5841) node[right]{$v_5(x_5,y_5)$} -- (0.3871,0.8381) node[right]{$v_3(x_3,y_3)$} -- (0.1754,0.9845) node[right,yshift=1ex]{$v_1(x_1,y_1)$} -- (0,1) node[above]{$v_{11}(0,1)$} -- (-0.1754,0.9845) node[left,yshift=1ex]{$v_{10}(x_{10},y_{10})$} -- (-0.3871,0.8381) node[left]{$v_8(x_8,y_8)$} -- (-0.5000,0.5851) node[left]{$v_6(x_6,y_6)$} -- (-0.4614,0.3088) node[left]{$v_4(x_4,y_4)$} -- (-0.2828,0.0957) node[left,yshift=-1ex]{$v_2(x_2,y_2)$} -- cycle;
		\draw (0,1) -- (0,0) -- (0.1754,0.9845) -- (-0.2828,0.0957) -- (0.3871,0.8381) -- (-0.4614,0.3088) -- (0.5000,0.5841) -- (-0.5000,0.5841) -- (0.4614,0.3088) -- (-0.3871,0.8381) -- (0.2828,0.0957) -- (-0.1754,0.9845) -- (0,0);
		\draw (0.0438,0.2461) arc (79.90:90.00:0.25) node[midway,above,xshift=0.2ex]{$\theta_0$};
		\draw (0.0608,0.7623) arc (242.73:259.90:0.25) node[midway,below]{$\theta_1$};
		\draw (-0.1154,0.2813) arc (47.94:62.73:0.25) node[midway,above,xshift=1.1ex]{$\theta_2$};
		\draw (0.1750,0.7058) arc (211.96:227.94:0.25) node[midway,below,xshift=-1.6ex]{$\theta_3$};
		\draw (-0.2210,0.3776) arc (15.98:31.96:0.25) node[midway,right,yshift=0.6ex]{$\theta_4$};
		\draw (0.2500,0.5841) arc (180.00:195.98:0.25) node[midway,left]{$\theta_5$};
	\end{tikzpicture}
\end{figure}

The \textit{skeleton} of a small polygon $P$ consists of the vertices of $P$, together with all of the line segments that connect two vertices of $P$ at unit distance from one another.
Let $n\geq6$ denote an even integer.
From \cite{Graham} it is known that the skeleton of an optimal $n$-gon $P$ forms a connected graph, and a \textit{linear thrackle}: each pair of line segments in the skeleton intersect one another, possibly at an endpoint.
Foster and Szabo \cite{FosterSzabo} proved that the skeleton of an optimal small polygon, considered as a graph, consists of an $(n-1)$-cycle, with a single additional pendant edge connected to the remaining vertex.
It is conjectured that this additional pendant edge forms an axis of symmetry in the skeleton of the optimal polygon; this is in fact the case for $n\leq12$.
We thus consider polygons having a skeleton of the form shown in Figure~\ref{figSkeleton}: a star with $n-1$ points on vertices $v_0$, \ldots, $v_{n-2}$, an additional vertex $v_{n-1}$ with distance $1$ from $v_0$, and the line connecting $v_0$ and $v_{n-1}$ forming an axis of symmetry for the polygon.


We place $v_0$ at the origin and $v_1$ at the point $(0,1)$ in $\mathbb{R}^2$.
Let $\theta_0$ denote the angle $\angle v_{n-1} v_0 v_1$, and for $0<k<n/2$ let $\theta_k=\angle v_{k-1} v_k v_{k+1}$.
Due to the symmetry of the construction, and the fact that the star forms a closed path, we have
\begin{equation}\label{eqnAngleSum}
\sum_{j=0}^{\frac{n}{2}-1} \theta_j = \frac{\pi}{2}.
\end{equation}
If we let $(x_k,y_k)$ denote the coordinates of $v_k$, then
\begin{equation}\label{eqnXY}
x_k = \sum_{j=0}^{k-1} (-1)^j \sin\biggl(\sum_{i=0}^j \theta_i\biggr), \quad
y_k = \sum_{j=0}^{k-1} (-1)^j \cos\biggl(\sum_{i=0}^j \theta_i\biggr).
\end{equation}
In addition, the vertices $v_{n/2-1}$ and $v_{n/2}$ are connected by a horizontal line in the skeleton, so
\begin{equation}\label{eqnXHalf}
x_{n/2-1} = -x_{n/2} = \frac{(-1)^{\frac{n}{2}}}{2}.
\end{equation}
We can compute the area $A=A(\theta_0,\ldots,\theta_{n/2-1})$ of such a polygon by determining the area of $n/2-1$ triangles $A_k$: $A_1$ is the area of the triangle $\Delta v_0 v_{n-1} v_1$, and $A_k$ is the area of $\Delta v_0 v_{k-1} v_{k+1}$ for $2\leq k<n/2$.
Then
\begin{equation}\label{eqnAreaDissection}
A = 2\sum_{k=1}^{\frac{n}{2}-1} A_k.
\end{equation}
It follows that
\begin{equation}\label{eqnGenlAk}
\begin{split}
2A_1 &= x_0 = \sin\theta_0,\\
2A_k &= x_{k+1}y_{k-1}-y_{k+1}x_{k-1}\\
&= \sin\theta_k + 2(-1)^k\left(x_k \sin\Bigl(\frac{\theta_k}{2}+\sum_{j=0}^{k-1} \theta_j\Bigr) + y_k\cos\Bigl(\frac{\theta_k}{2}+\sum_{j=0}^{k-1} \theta_j\Bigr)\right)\sin \frac{\theta_k}{2}\\
&= \sum_{i=0}^{k-2} (-1)^{i} \left(\sin \Bigl( \sum_{j=0}^{i+1} \theta_{k-j} \Bigr) - \sin \Bigl(\sum_{j=1}^{i+1} \theta_{k-j}\Bigr)\right)
\end{split}
\end{equation}
for $2\leq k<n/2$.
We thus obtain an expression for the area in terms of the $n/2$ angles $\theta_k$.

In \cite{BinganeA,Mossinghoff05}, this formulation was simplified to employ just three variables, $\alpha$, $\beta$, and $\gamma$, by taking
\begin{equation}\label{eqnAngles1}
\begin{split}
\theta_0 &= \alpha,\\
\theta_1 &= \beta+\gamma,\\
\theta_2 &= \beta-\gamma,\\
\theta_k &= \beta \quad \textrm{for $k\geq3$}.
\end{split}
\end{equation}
This configuration was selected to mimic some of the pattern observed in \cite{Mossinghoff05} for the small $n$-gons with large area for $n\leq20$, which were constructed using heuristic optimization methods over the parameters $\theta_0$, \ldots, $\theta_{n/2-1}$.
We extend these calculations in Section~\ref{secSmallN} for $n\leq120$ and note that this pattern continues: see Table~\ref{tableAngles}.
There, in each polygon constructed, the angles $\theta_i$ show a pattern of damped oscillation, with the odd-indexed values for the constructed $n$-gon appearing to converge from above to a limiting value in $(\frac{\pi}{n}, \frac{\pi}{n-1})$, and the even-indexed angles (after $\theta_0$) converging to the same value from below.
Using $\alpha$, $\beta$, and $\gamma$ in this way then allowed approximating the largest variations one appears to expect in the sequence of angles $\theta_i$, while keeping the analysis tractable by using only three variables.

As in \cite{BinganeA}, we note certain constraints on $\alpha$, $\beta$, and $\gamma$ inherited by the geometry, namely
\begin{gather}\label{eqnABRel}
\alpha+\left(\frac{n}{2}-1\right)\beta = \frac{\pi}{2},\\\label{eqnABCRel}
\sin(\alpha+\beta+\gamma) = \sin\alpha + \frac{\sin(\alpha+3\beta/2)}{2\cos(\beta/2)}.
\end{gather}
The former clearly follows from \eqref{eqnAngleSum}, and the latter is a consequence of this, combined with \eqref{eqnXY} and \eqref{eqnXHalf}.
After using these to eliminate $\beta$ and $\gamma$, the expression for the area in \cite{BinganeA} thus relied only on $\alpha$, and an asymptotic analysis was performed on this expression after setting $\alpha=a\pi/n+b\pi/n^2+c\pi/n^3$.
A similar analysis was performed in \cite{Mossinghoff05}.
Both of those works found that
\begin{equation}\label{eqnBoundBn}
\overline{A}_n - A(P_n) = \frac{(5303-456\sqrt{114})\pi^3}{5808n^3} + \frac{(192107-17934\sqrt{114})\pi^3}{21296n^4} + O\left(\frac{1}{n^5}\right),
\end{equation}
with $P_n=M_n$ or $B_n$ respectively, with
the analysis in \cite{BinganeA} producing an improvement in the $1/n^5$ term.

\section{Proof of Theorem~\ref{thmArea}}\label{secProof}

We obtain improved small polygons by generalizing the construction of Section~\ref{secPrior}, keeping some additional variables to allow for more variation in the sequence of angles $\theta_i$, beyond what is captured by \eqref{eqnAngles1}.
Let $r$ denote a positive integer, and suppose $n$ is even and $n \geq 2r+4$.
We describe a construction for a small $n$-gon which involves $r+2$ variables.

Assume first that $r$ is even.
Our variables are $\alpha$, $\beta$, $\beta_1$, \ldots, $\beta_{r/2}$, and $\gamma_1$, \ldots, $\gamma_{r/2}$, and we set
\begin{equation*}\label{eqnAnglesEven}
\begin{split}
\theta_0 &= \alpha,\\
\theta_{2i-1} &= \beta_i + \gamma_i,\\
\theta_{2i} &= \beta_i - \gamma_i, \hspace{1ex} \textrm{for $1\leq i\leq r/2$},\\
\theta_k &= \beta, \hspace{1ex} \textrm{for $r<k<n/2$}.
\end{split}
\end{equation*}
For convenience, let
\[
\varphi_r = \alpha + 2\sum_{i=1}^{r/2} \beta_i.
\]

We derive an expression for the area of the small $n$-gon in terms of $\alpha$, $\beta$, and the $\beta_i$ and $\gamma_i$.
For $k > r$, the coordinates $(x_k,y_k)$ in~\eqref{eqnXY} become
\begin{equation}\label{eqnXYModel}
\begin{split}
x_k
&= x_r + \sum_{j=r}^{k-1} (-1)^j \sin (\varphi_r + (j-r)\beta)\\
&= x_r + \frac{\sin \left(\varphi_r - \frac{\beta}{2}\right) - (-1)^k\sin \left(\varphi_r + (2k-2r-1)\frac{\beta}{2} \right)}{2\cos(\beta/2)},\\
y_k
&= y_r + \sum_{j=r}^{k-1} (-1)^j \cos (\varphi_r + (j-r)\beta)\\
&= y_r + \frac{\cos \left(\varphi_r - \frac{\beta}{2}\right) - (-1)^k\cos \left(\varphi_r + (2k-2r-1)\frac{\beta}{2} \right)}{2\cos(\beta/2)}
\end{split}
\end{equation}
when $r$ is even.
The constraint \eqref{eqnABRel} on the sum of the angles is now
\begin{equation}\label{eqnTotalAngle}
\varphi_r + \left(\frac{n}{2} - r-1\right)\beta = \frac{\pi}{2},
\end{equation}
and by combining this with \eqref{eqnXHalf} and \eqref{eqnXYModel}, we deduce a generalization of \eqref{eqnABCRel}:
\begin{equation}\label{eqnABCRelGenl}
x_r + \frac{\sin \left(\varphi_r - \frac{\beta}{2}\right)}{2\cos(\beta/2)} = 0.
\end{equation}
Using \eqref{eqnGenlAk} and \eqref{eqnXYModel}, the area $2A_k$ is then
\begin{align*}
	2A_k
	&= \sin \beta + (-1)^k 2\sin(\beta/2) \left(x_k\sin\left(\varphi_r+(2k-2r-1)\frac{\beta}{2}\right) \right.\\
	&\qquad\left. {}+ y_k \cos \left(\varphi_r+(2k-2r-1)\frac{\beta}{2}\right) \right)\\
	&= \sin \beta - \tan(\beta/2) + (-1)^k  2\sin(\beta/2)\left( x_r\sin\left(\varphi_r+(2k-2r-1)\frac{\beta}{2}\right) \right.\\
	&\qquad\left. {}+ y_r \cos \left(\varphi_r+(2k-2r-1)\frac{\beta}{2}\right) + \frac{\cos ((k-r)\beta)}{2\cos(\beta/2)}\right)
\end{align*}
for $r < k < n/2$.
Combining this with \eqref{eqnTotalAngle} and \eqref{eqnABCRelGenl}, it follows that
\[
\sum_{k=r+1}^{\frac{n}{2}-1} 2A_k = \left(\frac{n}{2}-r-1\right)\left(\sin \beta - \tan(\beta/2)\right) - \left(x_r\sin \varphi_r + y_r\cos \varphi_r + \frac{1}{2}\right) \tan(\beta/2).
\]
We thus obtain an expression for the area of our small $n$-gon having $r+2$ variables, whose number of terms depends only on $r$:
\begin{equation}\label{eqnArea}
A = \sum_{k=1}^r 2A_k + \left(\frac{n}{2}-r-1\right)\left(\sin \beta - \tan(\beta/2)\right)
 - \left(x_r\sin \varphi_r + y_r\cos \varphi_r + \frac{1}{2}\right) \tan(\beta/2).
\end{equation}
We then eliminate $\beta$ using \eqref{eqnTotalAngle} and $\gamma_{r/2}$ with \eqref{eqnABCRelGenl}.

If $r$ is odd, we employ the same strategy as with $r+1$, except we set $\beta_{\frac{r+1}{2}}=\beta$.
The scheme \eqref{eqnAngles1} therefore corresponds to the case $r=1$.

For each even $n\geq6$ and $r\geq1$, let $Q_{n,r}$ denote the small polygon obtained by maximizing our area function \eqref{eqnArea} over the parameters $\alpha$, $\beta_1$, \ldots, $\beta_{\floor{r/2}}$, $\gamma_1$, \ldots, $\gamma_{\ceiling{r/2}-1}$, for $\alpha\in[\frac{\pi}{2n-2},\frac{\pi}{n}]$, $\beta_i\in[\frac{\pi}{n},\frac{2\pi}{n}]$ for each $i$, and $\gamma_i\in[0,\frac{\pi}{n}]$ for each $i$.
An asymptotic analysis reveals that the area is maximized for large $n$ when
\begin{align*}
\alpha(n) &= \frac{a\pi}{n} + O\left(\frac{1}{n^2}\right),\\
\beta_i(n) &= \frac{b_i\pi}{n} + O\left(\frac{1}{n^2}\right), \quad \textrm{for $1\leq i\leq \floor{r/2}$},\\
\gamma_i(n) &= \frac{c_i\pi}{n} + O\left(\frac{1}{n^2}\right), \quad \textrm{for $1\leq i\leq \ceiling{r/2}-1$},
\end{align*}
where $a$, $b_1$, \ldots, $b_{\floor{r/2}}$, $c_1$, \ldots, $c_{\ceiling{r/2}-1}$ maximize a particular cubic polynomial in $r$ variables.
When $r=1$ this polynomial is
\[
-\frac{\pi^3}{192}(88 a^3 + 84 a^2  - 222 a +  107),
\]
while at $r=2$ it is
\[
-\frac{\pi^3}{192}(88 a^3 + 12a^2(8b_1-1) - 6a(16b_1^2+21)
 + 128 b_1^3 - 48 b_1^2  - 216 b_1 + 243),
\]
and $r=3$ yields
\[
\begin{aligned}
-&\frac{\pi^3}{192}( 88 a^3 + 12a^2(16 b_1 - 12 c_1 + 7) - 6a(32 b_1^2 + 64 b_1 c_1 - 80 c_1^2 + 56 c_1 + 37) \\
& + 128 b_1^3 + 192 b_1^2 c_1 + 384 b_1 c_1^2 - 384 c_1^3 +
336 c_1^2 - 240 b_1 + 204 c_1 + 267).
\end{aligned}
\]
After removing the factor $-\pi^3$, we use the \texttt{NMinimize} function in Mathematica to determine the optimal value for these polynomials by numerical methods, requiring $0 \le a \le 1$, $0 \le b_i \le 2$ for each $i$, and $0 \le c_i \le 1/3$ for each $i$.
This produces an asymptotic estimate for $A(Q_{n,r})$ having the form 
\begin{equation}\label{eqnAsympQ}
A(Q_{n,r}) = \frac{\pi}{4} - \frac{5\pi^3}{48n^2} - \frac{q_r \pi^3}{n^3} + O\left(\frac{1}{n^4}\right).
\end{equation}

For completeness we let $Q_{n,0}$ denote the polygon created by selecting $\alpha$ optimally in the $n$-gon with $\theta_0=\alpha$ and $\theta_i=\beta$ for $1\leq i<n/2$, subject to \eqref{eqnABRel}, so when $\alpha=\pi/(2n-2)$ and $\beta=\pi/(n-1)$.
This is the polygon created by simply adding a vertex at unit distance antipodal to one vertex of the regular small $(n-1)$-gon, as in Figure~\ref{figSkeleton} for $n=12$.
The polygon $Q_{n,0}$ then has
\[
q_0 = \frac{7}{48} = 0.1458333\ldots\,.
\]

For $r=1$, as reported in \cite{BinganeA,Mossinghoff05}, the optimal choice for $a$ is
\[
a = \frac{2\sqrt{114}-7}{22} = 0.652461659\ldots\,,
\]
which produces
\[
q_1 = \frac{5545-456\sqrt{114}}{5808} = 0.116434627\ldots\,.
\]
At $r=2$, we obtain
\[
q_2 = 0.115697150\ldots\,.
\]
In fact, $q_2$ is a root of the polynomial
\[
x^4 - \frac{70705}{15876}x^3 + \frac{269167127}{41150592}x^2 - \frac{3381027871}{987614208}x + \frac{737985313}{2341011456}.
\]
The case $r=3$ yields a further improvement,
\[
q_3 = 0.115089913\ldots\,,
\]
which is a root of a polynomial with rational coefficients and degree $8$:
\begin{equation*}
\tiny\begin{split}
x^8 &- \frac{338067189760423194}{232662255261540774} x^7 + \frac{1980606171874180754147}{22335576505107914304} x^6
 - \frac{158140620301705167575191}{536053836122589943296} x^5\\
 &+ \frac{59647522303796634759434731}{102922336535537269112832} x^4
 - \frac{836103610314364495378933003}{1235068038426447229353984} x^3 + \frac{52675103710698128327456883067}{118566531688938934017982464} x^2\\
 &- \frac{14538141342029184829034957803}{105392472612390163571539968} x + \frac{442235633612728385344035304147}{40470709483157822811471347712}.
\end{split}
\end{equation*}
We improve this further with subsequent values of $r$: our results through $r=16$ are summarized in Table~\ref{tableResults}.

The polygons $Q_{n,r}$ are small by construction, and for even $n\geq6$ we set
\[
Q_n = \begin{cases}
Q_{n,n/2-2} & n \leq 34,\\
Q_{n,16} & n \geq 36.
\end{cases}
\]
The first statement of Theorem~\ref{thmArea} then follows by combining \eqref{eqnAsympQ} at $r=16$ with \eqref{eqnAreaBound} and \eqref{eqnBoundBn}.
For the final statement, we calculate the areas of $Q_{6,1}$, $Q_{8,2}$, $Q_{10,3}$, and $Q_{12,4}$ by optimizing over $\alpha$ and the relevant $\beta_i$ and $\gamma_i$.
The values we obtain are consistent with the optimal areas for $n=6$ from \cite{Graham}, $n=8$ from \cite{AHMX}, and $n=10$ and $12$  from \cite{HenrionMessine}, and the polygons $Q_n$ are optimal in these cases.
Values for the area and angles of these polygons are recorded in Table~\ref{tableOptimalSmallN}.
\qed

\vskip\baselineskip

Additional small improvements can certainly be obtained using larger values for $r$.
Of course, the procedure becomes more computationally onerous as $r$ increases, due to the increasing complexity of the optimization procedure as the number of variables grows.
Such improvements are likely to be minuscule, however, given the rapid convergence exhibited in the $q_r$ values we compute.

\begin{table}[tbhp]
\caption{Optimal values of $q_r$ in \eqref{eqnAsympQ} for $Q_{n,r}$, together with values for the free parameters that produce this coefficient.}\label{tableResults}
\resizebox{\linewidth}{!}{
\begin{tabular}[c]{|r|ccL{1in}L{1.2in}|}\toprule
\TS\BS$r$ & $q_r$ & $a$ & $b_1,\ldots, b_{\floor{r/2}}$ & $c_1,\ldots, c_{\ceiling{r/2}-1}$\\\midrule
\small $0$ & \small $0.1458333333333333$ & & &\\\hline
\small$1$ & \small$0.1164346275953378$ & \small$0.6524616592755737$ & \tiny & \tiny\\\hline
\small$2$ & \small$0.1156971503834968$ & \small$0.6554858160170336$ & \tiny$1.022718374818576$ & \tiny\\\hline
\small$3$ & \small$0.1150899130453658$ & \small$0.6585214692355722$ & \tiny$1.027063969740726$ & \tiny$0.06794672543480737$\\\hline
\small$4$ & \small$0.1150687309140004$ & \small$0.6586249743177490$ & \tiny$1.027209190884176$, $1.003828109549754$ & \tiny$0.06761473542307868$\\\hline
\small$5$ & \small$0.1150557470337394$ & \small$0.6586900229138448$ & \tiny$1.027300358228207$, $1.004461819824590$ & \tiny$0.06740615925761905$, $0.01028318684355288$\\\hline
\small$6$ & \small$0.1150552764425733$ & \small$0.6586923731079715$ & \tiny$1.027303650624826$, $1.004484647927062$, $1.000570284354183$ & \tiny$0.06739862447575472$, $0.01023212957791077$\\\hline
\small$7$ & \small$0.1150549998390641$ & \small$0.6586937593365635$ & \tiny$1.027305592742427$, $1.004498111789651$, $1.000662623065463$ & \tiny$0.06739417976638237$, $0.01020201310839604$, $0.001509019841357918$\\\hline
\small$8$ & \small$0.1150549897593822$ & \small$0.6586938098307351$ & \tiny$1.027305663478150$, $1.004498602162219$, $1.000665984979524$, $1.000083456862159$ & \tiny$0.06739401787457336$, $0.01020091616496605$, $0.001501502526879076$\\\hline
\small$9$ & \small$0.1150549838708233$ & \small$0.6586938392297916$ & \tiny$1.027305704817829$, $1.004498888799399$, $1.000667949966222$, $1.000096926042045$ & \tiny$0.06739392319318871$, $0.01020027499791677$, $0.001497108650978320$, $0.0002203516777390261$\\\hline
\small$10$ & \small$0.1150549836560699$ & \small$0.6586938403067916$ & \tiny$1.027305706321311$, $1.004498899243624$, $1.000668021641018$, $1.000097417195005$, $1.000012181578676$ & \tiny$0.06739391975105442$, $0.01020025160734295$, $0.001496948418422726$, $0.0002192534425102402$\\\hline
\small$11$ & \small$0.1150549835307224$ & \small$0.6586938408223755$ & \tiny$1.027305707189345$, $1.004498905348772$, $1.000668063456477$, $1.000097703894909$, $1.000014146646888$ & \tiny$0.06739391767843365$, $0.01020023796370542$, $0.001496854886068534$, $0.0002186123676843603$, $0.00003215289654154845$\\\hline
\small$12$ & \small$0.1150549835261505$ & \small$0.6586938408430183$ & \tiny$1.027305707214692$, $1.004498905576070$, $1.000668064981190$, $1.000097714354882$, $1.000014218316758$, $1.000001777358190$ & \tiny$0.06739391761102890$, $0.01020023745865263$, $0.001496851472990505$, $0.0002185889876955666$, $0.00003199263572187769$\\\hline
\small$13$ & \small$0.1150549835234823$ & \small$0.6586938407444091$ & \tiny$1.027305707207800$, $1.004498905698189$, $1.000668065875004$, $1.000097720453880$, $1.000014260156537$, $1.000002064060629$ & \tiny$0.06739391752000037$, $0.01020023715826531$, $0.001496849480792517$, $0.0002185753423198028$, $0.00003189910564581757$, $4.691124054714891\cdot10^{-6}$\\\hline
\small$14$ & \small$0.1150549835233850$ & \small$0.6586938407702582$ & \tiny$1.027305707207799$, $1.004498905739078$, $1.000668065815119$, $1.000097720770217$, $1.000014261621179$, $1.000002074524450$, $1.000000259301370$ & \tiny$0.06739391752829402$, $0.01020023716420131$, $0.001496849406827393$, $0.0002185748364450190$, $0.00003189570671658796$, $4.667741696766253\cdot10^{-6}$\\\hline
\small$15$ & \small$0.1150549835233282$ & \small$0.6586938406842894$ & \tiny$1.027305707182444$, $1.004498905675305$, $1.000668065987254$, $1.000097720721770$, $1.000014262600232$, $1.000002080718492$, $1.000000301103361$ & \tiny$0.06739391752641309$, $0.01020023709721798$, $0.001496849396193236$, $0.0002185745506217027$, $0.00003189368791732473$, $4.654108249822218\cdot10^{-6}$, $6.844384307984062\cdot10^{-7}$\\\hline
\small$16$ & \small$0.1150549835233261$ & \small$0.6586938406334803$ & \tiny$1.027305707190814$, $1.004498905736593$, $1.000668065901726$, $1.000097720834017$, $1.000014262573901$, $1.000002080858254$, $1.000000302668948$, $1.000000037787016$ & \tiny$0.06739391746458313$, $0.01020023714027207$, $0.001496849368106921$, $0.0002185745455879002$, $0.00003189363350580009$, $4.653599949365328\cdot10^{-6}$, $6.810134779486807\cdot10^{-7}$\\\bottomrule
\end{tabular}
}
\end{table}

\begin{table}[tbp]
\caption{Constructing optimal polygons for small $n$.}\label{tableOptimalSmallN}
\footnotesize
\centering
\resizebox{\linewidth}{!}{%
\begin{tabular}[c]{|c|ccll|}\toprule
\TS\BS$n$ & $A(Q_{n,n/2-2})$ & $\alpha$ & \multicolumn{1}{c}{$\beta_1$, $\beta_2$} & \multicolumn{1}{c|}{$\gamma_1$}\\\midrule
6 & $0.6749814429301047$ & $0.3509301888703616$ &&\\
8 & $0.7268684827516268$ & $0.2652408674910718$ & $0.4379295350493946$ &\\
10 & $0.7491373458778303$ & $0.2126101953284637$ & $0.3433714044229845$ & $0.02476000789351616$\\
12 & $0.7607298734487962$ & $0.1770854623284314$ & $0.2827755557037131$, & $0.01982894085863103$\\
& & & $0.2763754214389234$ &\\\bottomrule
\end{tabular}%
}
\end{table}

\section{Constructions for small $n$}\label{secSmallN}

We construct some small polygons with large area for particular values of $n$ and display their values in two tables.
First, for even integers $n$ with $6\leq n\leq 120$, we calculate the small $n$-gon with maximal area $P_n^*$, assuming the presence of an axis of symmetry.
We employ the skeleton of Foster and Szabo as in Figure~\ref{figSkeleton}, and assume that the polygon is symmetric about the line connecting $v_0$ and $v_{n-1}$.
For each such $n$, using \eqref{eqnAreaDissection} and \eqref{eqnGenlAk} we construct $P_n^*$ by maximizing the area $A$ over $n/2$ variables $\theta_0, \theta_1, \ldots, \theta_{\frac{n}{2}-1}$, subject to \eqref{eqnAngleSum} and~\eqref{eqnXHalf}.
More precisely,
\begin{equation}\label{eqnMaximalAreas}
\begin{split}
A(P_n^*) = \max_{\theta_0, \theta_1, \ldots, \theta_{\frac{n}{2}-1}} \quad
& \sin \theta_0 + \sum_{k=2}^{\frac{n}{2}-1} 2A_k(\theta_1,\theta_2, \ldots, \theta_k)\\
\subj \quad
& \sum_{k=0}^{\frac{n}{2}-1}\theta_k = \frac{\pi}{2},\\
&\sum_{i=0}^{\frac{n}{2}-2} (-1)^{i} \sin \Bigl(\sum_{j=0}^{i}\theta_j\Bigr) = \frac{(-1)^{\frac{n}{2}}}{2},\\
& 0 \le \theta_0 \le \pi/6,\\
& 0 \le \theta_k \le \pi/3, \quad 1\leq k \leq n/2-1.
\end{split}
\end{equation}

Problem~\eqref{eqnMaximalAreas} was solved on the NEOS Server~6.0 using AMPL with the nonlinear programming solver Ipopt~3.13.4 \cite{Ipopt}.
The AMPL code is available in OPTIGON~\cite{Optigon}, a free package for extremal convex small polygons available on GitHub.
For selected even $n\leq120$, Table~\ref{tableAngles} shows the values $\theta_i^*$ that we calculated for constructing $P_n^*$,
and each value $A(P_n^*)$ is displayed in Table~\ref{tableSmallN}.
The areas shown here for $n\leq20$ agree with those from \cite{Mossinghoff05}, and the values of $A(P_n^*)$ for larger $n$ in Table~\ref{tableSmallN} match or slightly exceed the best value found in the literature \cite{BinganeLSP,Pinter,PKC}.
Ipopt required less than $1$ second to compute each value in Table~\ref{tableSmallN}.

Second, for selected even $n\leq120$ we determine the area of $Q_{n,r}$ for $0\leq r\leq4$ by optimizing \eqref{eqnArea} over the $r$ parameters $\alpha$, $\beta_1$, \ldots, $\beta_{\floor{r/2}}$, $\gamma_1$, \ldots, $\gamma_{\ceiling{r/2}-1}$.
These areas are also displayed in Table~\ref{tableSmallN}, along with that of the regular $n$-gon $R_n$ and the upper bound $\overline{A}_n$.
Julia and MATLAB functions that give the coordinates of the vertices of all polygons presented in this work are provided in OPTIGON\@.

In all tables in this paper, each numerical value is rounded at the last displayed digit.

\begin{table}[tbhp]
	\footnotesize
	\centering
	\caption{Angles $\theta_0^*$, $\theta_1^*$, \ldots, $\theta_{\frac{n}{2}-1}^*$ of $P_n^*$.}\label{tableAngles}
	\resizebox{\linewidth}{!}{
		\begin{tabular}{|c|ccccccccc|}
			\toprule
			$n$ & $i$ & $\theta_{8i}^*$ & $\theta_{8i+1}^*$ & $\theta_{8i+2}^*$ & $\theta_{8i+3}^*$ & $\theta_{8i+4}^*$ & $\theta_{8i+5}^*$ & $\theta_{8i+6}^*$ & $\theta_{8i+7}^*$ \\
			\midrule
			6	&	0	&	0.350930	&	0.653342	&	0.566524	&&&&&\\\hline
			8	&	0	&	0.265241	&	0.470631	&	0.405228	&	0.429696	&&&&\\\hline
			10	&	0	&	0.212610	&	0.368131	&	0.318611	&	0.339137	&	0.332306	&&&\\\hline
			12	&	0	&	0.177085	&	0.302604	&	0.262947	&	0.279461	&	0.273290	&	0.275409	&&\\\hline
			14	&	0	&	0.151583	&	0.257026	&	0.233904	&	0.237628	&	0.232444	&	0.234442	&	0.233769	&\\\hline
			16	&	0	&	0.132428	&	0.223448	&	0.194967	&	0.206716	&	0.202285	&	0.204013	&	0.203359	&	0.203580	\\\hline
			18	&	0	&	0.117533	&	0.197661	&	0.172654	&	0.182938	&	0.179070	&	0.180577	&	0.179999	&	0.180218	\\
			&	1	&	0.180145	&&&&&&&\\\hline
			20	&	0	&	0.105629	&	0.177228	&	0.154925	&	0.164076	&	0.160640	&	0.161977	&	0.161464	&	0.161661	\\
			&	1	&	0.161586	&	0.161611	&&&&&&\\\hline
			22	&	0	&	0.0959016	&	0.160633	&	0.140496	&	0.148745	&	0.145652	&	0.146854	&	0.146393	&	0.146571	\\
			&	1	&	0.146503	&	0.146528	&	0.146520	&&&&&\\\hline
			24	&	0	&	0.0878067	&	0.146886	&	0.128526	&	0.136037	&	0.133224	&	0.134316	&	0.133898	&	0.134059	\\
			&	1	&	0.133997	&	0.134021	&	0.134012	&	0.134015	&&&&\\\hline
			30	&	0	&	0.0700443	&	0.116891	&	0.102361	&	0.108291	&	0.106075	&	0.106933	&	0.106605	&	0.106731	\\
			&	1	&	0.106683	&	0.106701	&	0.106694	&	0.106697	&	0.106696	&	0.106696	&	0.106696	&\\\hline
			40	&	0	&	0.0523626	&	0.0872236	&	0.0764267	&	0.0808253	&	0.0791841	&	0.0798182	&	0.0795763	&	0.0796689	\\
			&	1	&	0.0796334	&	0.0796470	&	0.0796418	&	0.0796437	&	0.0796429	&	0.0796432	&	0.0796431	&	0.0796431	\\
			&	2	&	0.0796431	&	0.0796431	&	0.0796431	&	0.0796431	&&&&\\\hline
			50	&	0	&	0.0418008	&	0.0695718	&	0.0609765	&	0.0644752	&	0.0631706	&	0.0636742	&	0.0634822	&	0.0635556	\\
			&	1	&	0.0635275	&	0.0635382	&	0.0635340	&	0.0635355	&	0.0635349	&	0.0635351	&	0.0635350	&	0.0635350	\\
			&	2	&	0.0635350	&	0.0635350	&	0.0635349	&	0.0635349	&	0.0635349	&	0.0635349	&	0.0635349	&	0.0635349	\\
			&	3	&	0.0635349	&&&&&&&\\\hline
			60	&	0	&	0.0347820	&	0.0578637	&	0.0507223	&	0.0536280	&	0.0525448	&	0.0529626	&	0.0528033	&	0.0528642	\\
			&	1	&	0.0528408	&	0.0528496	&	0.0528462	&	0.0528474	&	0.0528469	&	0.0528470	&	0.0528469	&	0.0528469	\\
			&	2	&	0.0528468	&	0.0528468	&	0.0528468	&	0.0528467	&	0.0528467	&	0.0528467	&	0.0528467	&	0.0528467	\\
			&	3	&	0.0528467	&	0.0528466	&	0.0528466	&	0.0528466	&	0.0528466	&	0.0528466	&&\\\hline
			70	&	0 &	0.0297804	&	0.0495296	&	0.0434206	&	0.0459055	&	0.0449793	&	0.0453364	&	0.0452002	&	0.0452521	\\
			& 1 &	0.0452321	&	0.0452396	&	0.0452366	&	0.0452376	&	0.0452371	&	0.0452372	&	0.0452371	&	0.0452370	\\
			& 2 &	0.0452370	&	0.0452369	&	0.0452369	&	0.0452369	&	0.0452368	&	0.0452368	&	0.0452368	&	0.0452367	\\
			& 3 &	0.0452367	&	0.0452367	&	0.0452367	&	0.0452367	&	0.0452366	&	0.0452366	&	0.0452366	&	0.0452366	\\
			& 4 &	0.0452366	&	0.0452366	&	0.0452366	&&&&&\\\hline
			80	&	0 &	0.0260361	&	0.0432947	&	0.0379568	&	0.0401276	&	0.0393185	&	0.0396302	&	0.0395112	&	0.0395565	\\
			& 1 &	0.0395390	&	0.0395454	&	0.0395428	&	0.0395437	&	0.0395432	&	0.0395432	&	0.0395431	&	0.0395430	\\
			& 2 &	0.0395430	&	0.0395429	&	0.0395429	&	0.0395428	&	0.0395428	&	0.0395427	&	0.0395427	&	0.0395426	\\
			& 3 &	0.0395426	&	0.0395426	&	0.0395425	&	0.0395425	&	0.0395425	&	0.0395424	&	0.0395424	&	0.0395424	\\
			& 4 &	0.0395424	&	0.0395424	&	0.0395424	&	0.0395424	&	0.0395424	&	0.0395423	&	0.0395423	&	0.0395423	\\\hline
			90	&	0 &	0.0231282	&	0.0384545	&	0.0337147	&	0.0356421	&	0.0349236	&	0.0352001	&	0.0350945	&	0.0351345	\\
			& 1 &	0.0351190 	&	0.0351247	&	0.0351223	&	0.0351230	&	0.0351226	&	0.0351225	&	0.0351224	&	0.0351223	\\
			& 2 &	0.0351222	&	0.0351221	&	0.0351221	&	0.0351220	&	0.0351219	&	0.0351219	&	0.0351218	&	0.0351218	\\
			& 3 &	0.0351217	&	0.0351217	&	0.0351216	&	0.0351216	&	0.0351215	&	0.0351215	&	0.0351215	&	0.0351214	\\
			& 4 &	0.0351214	&	0.0351214	&	0.0351214	&	0.0351213	&	0.0351213	&	0.0351213	&	0.0351213	&	0.0351213	\\
			& 5 &	0.0351213	&	0.0351213	&	0.0351213	&	0.0351213	&	0.0351213	&&&\\\hline
			100	&	0 &	0.0208046	&	0.0345883	&	0.0303258	&	0.0320589	&	0.0314127	&	0.0316612	&	0.0315662	&	0.0316021	\\
			& 1 &	0.0315881	&	0.0315931	&	0.0315909	&	0.0315915	&	0.0315911	&	0.0315910	&	0.0315909	&	0.0315907	\\
			& 2 &	0.0315906	&	0.0315905	&	0.0315904	&	0.0315903	&	0.0315903	&	0.0315902	&	0.0315901	&	0.0315900	\\
			& 3 &	0.0315900	&	0.0315899	&	0.0315898	&	0.0315898	&	0.0315897	&	0.0315897	&	0.0315897	&	0.0315896	\\
			& 4 &	0.0315896	&	0.0315895	&	0.0315895	&	0.0315895	&	0.0315894	&	0.0315894	&	0.0315894	&	0.0315894	\\
			& 5 &	0.0315893	&	0.0315893	&	0.0315893	&	0.0315893	&	0.0315893	&	0.0315893	&	0.0315893	&	0.0315893	\\
			& 6 &	0.0315893	&	0.0315893	&&&&&&\\\hline
			110	&	0 &	0.0189055	&	0.0314290	&	0.0275563	&	0.0291308	&	0.0285436	&	0.0287692	&	0.0286828	&	0.0287153	\\
			& 1 &	0.0287025	&	0.0287070	&	0.0287050	&	0.0287054	&	0.0287050	&	0.0287049	&	0.0287048	&	0.0287046	\\
			& 2 &	0.0287045	&	0.0287043	&	0.0287042	&	0.0287041	&	0.0287040	&	0.0287039	&	0.0287038	&	0.0287037	\\
			& 3 &	0.0287037	&	0.0287036	&	0.0287035	&	0.0287034	&	0.0287034	&	0.0287033	&	0.0287033	&	0.0287032	\\
			& 4 &	0.0287031	&	0.0287031	&	0.0287031	&	0.0287030	&	0.0287030	&	0.0287029	&	0.0287029	&	0.0287029	\\
			& 5 &	0.0287028	&	0.0287028	&	0.0287028	&	0.0287028	&	0.0287027	&	0.0287027	&	0.0287027	&	0.0287027	\\
			& 6 &	0.0287027	&	0.0287027	&	0.0287026	&	0.0287026	&	0.0287026	&	0.0287026	&	0.0287026	&\\\hline
			120	&	0 &	0.0173243	&	0.0287991	&	0.0252508	&	0.0266933	&	0.0261551	&	0.0263616	&	0.0262824	&	0.0263121	\\
			& 1 &	0.0263003	&	0.0263043	&	0.0263024	&	0.0263028	&	0.0263023	&	0.0263022	&	0.0263020	&	0.0263018	\\
			& 2 &	0.0263017	&	0.0263015	&	0.0263014	&	0.0263012	&	0.0263011	&	0.0263010	&	0.0263009	&	0.0263008	\\
			& 3 &	0.0263007	&	0.0263006	&	0.0263005	&	0.0263004	&	0.0263003	&	0.0263003	&	0.0263002	&	0.0263001	\\
			& 4 &	0.0263001	&	0.0263000	&	0.0262999	&	0.0262999	&	0.0262998	&	0.0262998	&	0.0262997	&	0.0262997	\\
			& 5 &	0.0262996	&	0.0262996	&	0.0262996	&	0.0262995	&	0.0262995	&	0.0262995	&	0.0262994	&	0.0262994	\\
			& 6 &	0.0262994	&	0.0262994	&	0.0262994	&	0.0262993	&	0.0262993	&	0.0262993	&	0.0262993	&	0.0262993	\\
			& 7 &	0.0262993	&	0.0262993	&	0.0262993	&	0.0262993	&&&&\\
			\bottomrule
		\end{tabular}
	}
\end{table}

\begin{table}[tbhp]
	\footnotesize
	\centering
	\caption{Comparing areas of small polygons.}\label{tableSmallN}
	\resizebox{\linewidth}{!}{
		\begin{tabular}{@{}c|cccccccc@{}}
			\toprule
			$n$ & $A(R_n)$ & $A(Q_{n,0})$ & $A(Q_{n,1})$ & $A(Q_{n,2})$ & $A(Q_{n,3})$ & $A(Q_{n,4})$ & $A(P_n^*)$ & $\overline{A}_n$ \\
			\midrule
			6	&	0.6495190528	&	0.6722882584	&	0.6749814429		&	--		&	--	&	--	&	0.6749814429	&	0.6877007594	\\
			8	&	0.7071067812	&	0.7253199909	&	0.7268542719	&	0.7268684828		&	--		&	--	&	0.7268684828	&	0.7318815691	\\
			10	&	0.7347315654	&	0.7482573378	&	0.7491189262	&	0.7491297887		&	0.7491373459		&	--	&	0.7491373459	&	0.7516135587	\\
			12	&	0.7500000000	&	0.7601970055	&	0.7607153082	&	0.7607228359		&	0.7607297471		&	0.7607298734	&	0.7607298734	&	0.7621336536	\\
			14	&	0.7592965435	&	0.7671877750	&	0.7675203660	&	0.7675256353		&	0.7675308404		&	0.7675309615	&	0.7675310111	&	0.7684036467	\\
			16	&	0.7653668647	&	0.7716285345	&	0.7718535572	&	0.7718573456		&	0.7718611688		&	0.7718612660	&	0.7718613220	&	0.7724408116	\\
			18	&	0.7695453225	&	0.7746235089	&	0.7747824059	&	0.7747852057		&	0.7747880405		&	0.7747881160	&	0.7747881651	&	0.7751926059	\\
			20	&	0.7725424859	&	0.7767382147	&	0.7768543958	&	0.7768565173		&	0.7768586570		&	0.7768587158	&	0.7768587560	&	0.7771522071	\\
			22	&	0.7747645313	&	0.7782865351	&	0.7783739622	&	0.7783756055		&	0.7783772514		&	0.7783772976	&	0.7783773302	&	0.7785970008	\\
			24	&	0.7764571353	&	0.7794540033	&	0.7795213955	&	0.7795226929		&	0.7795239821		&	0.7795240189	&	0.7795240452	&	0.7796927566	\\
			30	&	0.7796688406	&	0.7816380102	&	0.7816725130	&	0.7816732130		&	0.7816738921		&	0.7816739122	&	0.7816739269	&	0.7817597927	\\
			40	&	0.7821723252	&	0.7833076096	&	0.7833221318	&	0.7833224422		&	0.7833227341		&	0.7833227431	&	0.7833227495	&	0.7833587784	\\
			50	&	0.7833327098	&	0.7840695435	&	0.7840769608	&	0.7840771244		&	0.7840772750		&	0.7840772797	&	0.7840772830	&	0.7840956746	\\
			60	&	0.7839634745	&	0.7844798073	&	0.7844840910	&	0.7844841875		&	0.7844842749		&	0.7844842777	&	0.7844842796	&	0.7844949027	\\
			70	&	0.7843439529	&	0.7847256986	&	0.7847283918	&	0.7847284534		&	0.7847285085		&	0.7847285103	&	0.7847285115	&	0.7847351925	\\
			80	&	0.7845909573	&	0.7848845934	&	0.7848863952	&	0.7848864368		&	0.7848864738		&	0.7848864750	&	0.7848864758	&	0.7848909473	\\
			90	&	0.7847603296	&	0.7849931681	&	0.7849944322	&	0.7849944617		&	0.7849944876		&	0.7849944885	&	0.7849944890	&	0.7849976272	\\
			100	&	0.7848814941	&	0.7850706272	&	0.7850715479	&	0.7850715695		&	0.7850715884		&	0.7850715890	&	0.7850715895	&	0.7850738759	\\
			110	&	0.7849711494	&	0.7851278167	&	0.7851285079	&	0.7851285242		&	0.7851285384		&	0.7851285389	&	0.7851285392	&	0.7851302562	\\
			120	&	0.7850393436	&	0.7851712379	&	0.7851717699	&	0.7851717826		&	0.7851717935		&	0.7851717939	&	0.7851717941	&	0.7851731162	\\
			\bottomrule
		\end{tabular}
	}
\end{table}

\bibliographystyle{amsplain}

\end{document}